\documentclass[a4paper,11pt]{article}
\usepackage{amscd}
\usepackage{amssymb}
\usepackage{amsthm}
\usepackage{xypic}
\usepackage{txfonts}
\usepackage{pstricks}

\let\O\relax
\usepackage{amsmath}
\advance\oddsidemargin-0.35in
\newcommand{\C}{\mathfrak C}
\newcommand{\F}{\varmathbb F}
\newcommand{\O}{\mathcal O}
\newcommand{\B}{\mathcal B}
\newcommand{\M}{\mathcal M}
\newcommand{\Spec}{\mbox{Spec}}
\newcommand{\Z}{\varmathbb Z}
\newcommand{\N}{\varmathbb N}
\newcommand{\R}{\varmathbb R}
\newcommand{\Q}{\varmathbb Q}
\newcommand{\Frac}{\mbox{Frac}}
\newcommand{\perf}{\mbox{perf}}
\renewcommand{\red}{\mbox{red}}

\begin{document}
\begin{center}\textbf{\Large Purity Results on $F$-crystals}\\
\vskip 0.5cm
Jinghao Li
\end{center}

\vskip0.5cm

Let $n\in \N$. Let $k$ be a perfect field of characteristic $p>0$. Let W(k) be the ring of Witt vectors with coefficient in $k$ and let $B(k)$ be its fractional field.\\

Let $\sigma$ be the absolute Frobenius automorphism  :\begin{center} $W(k)\xrightarrow{\sigma} W(k)$\\

$x=(x_0,x_1,\cdots)\mapsto \sigma (x)=(x_0^p,x_1^p,\cdots)$.\\

\end{center}

We extend $\sigma$ naturally to an automorphism of $B(k)$.\\

A $\sigma^n$-$F$-crystal (or $F^n$-crystal or $F-$crystal if $n=1$) over $k$ (or Spec $k$) is a pair $(M,F)$ consisting of a free $W(k)$-module M of finite rank, together with a $\sigma^n$-linear endomorphism $F:M\rightarrow M$, i.e. $F$ is additive and $F(\lambda x)=\sigma^n (\lambda)F(x)~ \forall \lambda\in W(k), \forall x\in M$, which induces an automorphism of $M\otimes_{\varmathbb{Z}_p}\varmathbb {Q}_p$.\\ 

The exterior powers of an $F^n$-crystal $(M,F)$ are the $F^{in}$-crystals $(\wedge^iM,\wedge^iF)~~(i=0,1,2,\cdots)$ with underlying modules $\wedge^iM$ and $\sigma^{in}$-endomorphisms $\wedge^i(F)$ defined by \\
$\wedge^i(F)(m_1\wedge\cdots\wedge m_i)=F(m_1)\wedge...\wedge F(m_i)$. For $i=0$ and $(M,F)\neq 0$, we define $(\wedge^0M,\wedge^0F)$ to be $(W(k),\sigma)$.\\

A morphism of $F^n$-crystals $f: (M,F)\rightarrow (M',F')$ is a $W(k)$-linear map $f:M\rightarrow M'$ such that $fF=F'f$.\\

The category of $F^n$-crystals up to isogeny is obtained from the category of $F^n$-crystals by keeping the same objects, but tensoring \textbf{Hom} groups, which are $\Z_p$ modules, over $\Z_p$ with $\Q_p$. An isogeny between $F^n$-crystals is a morphism of $F^n$-crystals which becomes an isomorphism in the category of $F^n$-crystals up to isogeny.\\

An $F^n$-crystal is said to be divisible by $\lambda>0$ if for all $m\in \N,$ we have $F^m=0$ mod $p^{[m\lambda]}$, where $[\cdot]:\R\to\Z$ is the floor function.\\

The Hodge slopes of an $F^n$-crystal $(M,F)$ are the integers defined as follows. The image $F(M)$ is a $W(k)$-submodule of $M$ of rank $r=\mbox{rank} (M)$, so by the theory of elementary divisors, since $W(k)$ is a discrete valuation ring, in particular a principal ideal domain, there exist $W(k)$-bases $\{v_1,...,v_r\}$ and $\{w_1,...,w_r\}$ of $M$ such that for all $1\le i\le r$ we have\\
$$F(v_i)=p^{a_i}w_i$$ for some integers $0\le a_1\le a_2\le\cdots\le a_r$. These integers are called the Hodge slopes of $(M,F)$.\\

The Hodge polygon of $(M,F)$ is the graph of the Hodge function on $[0,r]$ defined on integers $0\le i \le r$ by\\
\begin{center} \[\mbox{Hodge}_F(i)=\mbox{Least Hodge slope of }(\wedge^iM,\lambda^iF)=\begin{cases} 0, & i=0,\\  a_1+\cdots+a_i, & 1\le i \le r\end{cases}\] \end{center}
and then extended linearly between successive integers.\\

The Newton slopes of an $F^n$-crystal $(M,F)$ are the sequence of $r=\mbox{rank}(M)$ rational numbers $0\le \lambda_1\le\cdots\le \lambda_r$ defined in the following way.\\

Pick an algebraically closed field extension $k'$ of $k$, and consider the $F^n$-crystal over $k'$: 
$(M\otimes_{W(k)}W(k'),F\otimes \sigma^n)$ obtained from $(M,F)$ by extension of scalars. For each 
non-negative rational number $\lambda=a/b$ with $a\in \Z,b\in \N$, $(a,b)=1$, we denote by $E(\lambda)$ the $F^n$-crystal over $k'$ defined by:\\

\begin{center} $E(\lambda)=((\varmathbb{Z}_p[T]/(T^b-p^a))\otimes_{\varmathbb{Z}_p}W(k')$, (multiplication by $T^n)\otimes \sigma^n$).\end{center}

\textbf {Theorem 0 (Dieudonn\'e--Manin, [De72, (Theorem, Page 85)])}
{\it If $k'$ is an algebraically closed extension of $k$, then for any $F^n$-crystal $(M,F)$ over $k$, $(M\otimes_{W(k)}B(k'),F\otimes \sigma^n)$ is isomorphic to a finite direct sum of $ E(\lambda_i)\otimes_{W(k')} B(k')$'s with $\lambda_i\in \Q_{\ge 0}$.}\\

By the above theorem, we can write $(M\otimes_{W(k)}B(k'),F\otimes \sigma^n)\cong \bigoplus\limits_{i=1}^s E(a_i/b_i)\otimes B(k')$ for an increasing sequence $a_1/b_1\le a_2/ b_2\le\cdots\le a_s/b_s$ with $\sum\limits_{i=1}^s b_i=r$. The Newton slopes of $(M,F)$ are defined to be the sequence of $r$ rational numbers $(\lambda_1, \cdots,\lambda_r):=(a_1/b_1$ repeated $b_1$ times, $a_2/b_2$ repeated $b_2$ times,$\cdots$, $a_s/b_s$ repeated $b_s$ times).\\

The Newton polygon of $(M,F)$ is the graph of the Newton function on $[0,r]$, defined on integers $0\le i \le r$ by

\begin{center} \[\mbox{Newton}_F(i)=\mbox{least Newton slope of }(\wedge^iM,\lambda^iF)=\begin{cases} 0, & i=0,\\  \lambda_1+\cdots+\lambda_i, & 1\le i \le r\end{cases}\] \end{center}

\noindent and then extended linearly between successive integers. Let $\{\mu_1,\mu_2,\cdots,\mu_t\}=\{\lambda_1,\lambda_2,\cdots,\lambda_r\}$ with $\mu_1<\mu_2<\cdots<\mu_t$ for $1\le t\le r$ and let $r_i$ be the multiplicity of $\mu_i$ for $1\le i\le t$.\\

Here is the graph of the Newton polygon of $(M,F)$:

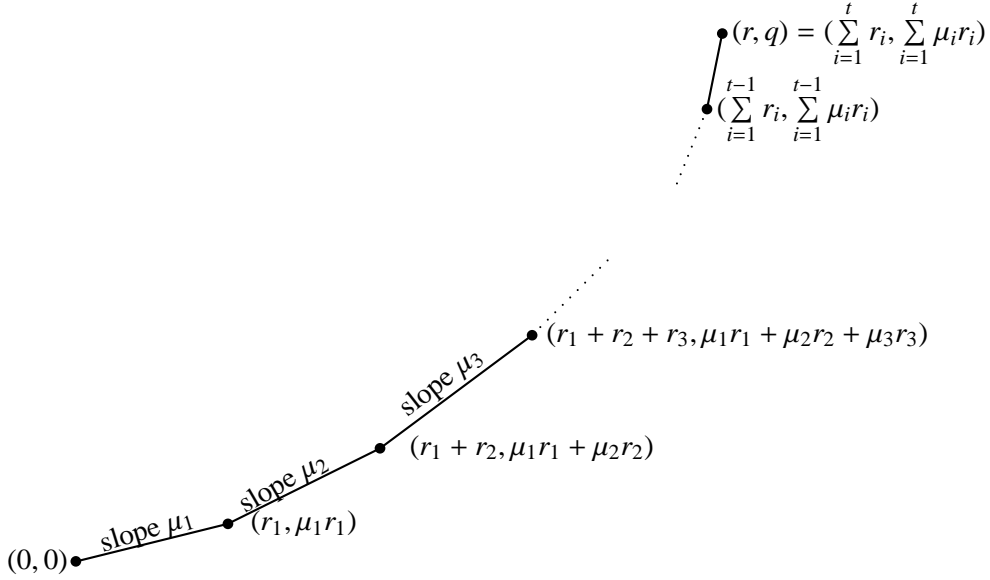
\begin{pspicture}(10,7)
\psdots[dotsize=4pt](0,0)(2,0.5)(4,1.5)(6,3)(8.5,7)(8.3,6)
\psline(0,0)(2,0.5)
\psline(2,0.5)(4,1.5) 
\psline(4,1.5)(6,3)
\psline[linestyle=dotted](6,3)(7,4)
\psline(8.3,6)(8.5,7)
\psline[linestyle=dotted](7.9,5)(8.3,6)
\rput[tr]{17}(1.5,0.75){\text{slope $\mu_1$}}
\rput[tr]{30}(3.2,1.5){\text{slope $\mu_2$}}
\rput[tr]{38}(5.2,2.9){\text{slope $\mu_3$}} 
\rput(-0.5,0){$(0,0)$}
\rput(3,0.5){$(r_1,\mu_1r_1)$}
\rput(6,1.5){$(r_1+r_2,\mu_1r_1+\mu_2r_2)$}
\rput(8.7,3){$(r_1+r_2+r_3,\mu_1r_1+\mu_2r_2+\mu_3r_3)$}
\rput(10.3,7){$(r,q)=(\sum\limits_{i=1}^t r_i,\sum\limits_{i=1}^t\mu_ir_i)$}
\rput(9.5,6){$(\sum\limits_{i=1}^{t-1} r_i,\sum\limits_{i=1}^{t-1}\mu_ir_i)$}
\end{pspicture} 

\vskip 1cm

For another characterization of the Newton slopes, we choose an auxiliary integer $N\ge 1$ which is divisible by $r!$, where $r=\mbox{rank}(M)$, and consider the discrete valuation ring $R=W(k')[X]/(X^N-p)=W(k')[p^{1/N}]$. We can extend $\sigma$ to an automorphism of $R$ by requiring that $\sigma (X+(X^N-p))=X+(X^N-p)$. For each $\lambda\in\frac{1}{N}\Z$, we can speak about $p^\lambda=X^{N\lambda}+(X^N-p)$ in $R$.\\

Let $K=\mbox{Frac}(R)$. By an analogue of Dieudonn\'e--Manin's theorem over $K$, we know that $M\otimes_{W(k)}K$ admits a $K$-basis $e_1,...,e_r$ which transforms under the $\sigma-$linear endomorphism $F\otimes \sigma$ by the formula
$(F\otimes \sigma)(e_i)=p^{\lambda_i}e_i$ for $1\le i\le r$. An equivalent characterization of the Newton slopes is by the existence of an $R$-basis $u_1,...,u_r$ of $M\otimes_{W(k)}R$ with respect to which the matrix of $F\otimes \sigma$ is upper-triangular, with $p^{\lambda_i}$'s along the diagonal, i.e. $F(u_i)\equiv p^{\lambda_i}u_i$ mod $\sum\limits_{1\le j<i} Ru_j$ for all $1\le i\le r$.\\

The second characterization of Newton slopes shows:\\

\noindent 1. The Newton slopes of the $m^{th}$ iterate $(M,F^m)$ of $(M,F)$ are $(m\lambda_1,...,m\lambda_r)$.\\
2. All Newton slopes $\lambda_i$ of $(M,F)$ are equal to $0$ if and only if $F$ is a $\sigma^n$-linear automorphism of $M$.\\
3. All Newton slopes $\lambda_i$ of $(M,F)$ are $>0$ if and only if $F$ is topologically nilpotent on $M$, i.e. if and only if we have $F^r(M)\subset pM$ where $r=\mbox{rank}(M)$.\\

The break points of an $F^n$-crystal $(M,F)$ are defined to be the points where the Newton polygon changes slopes.

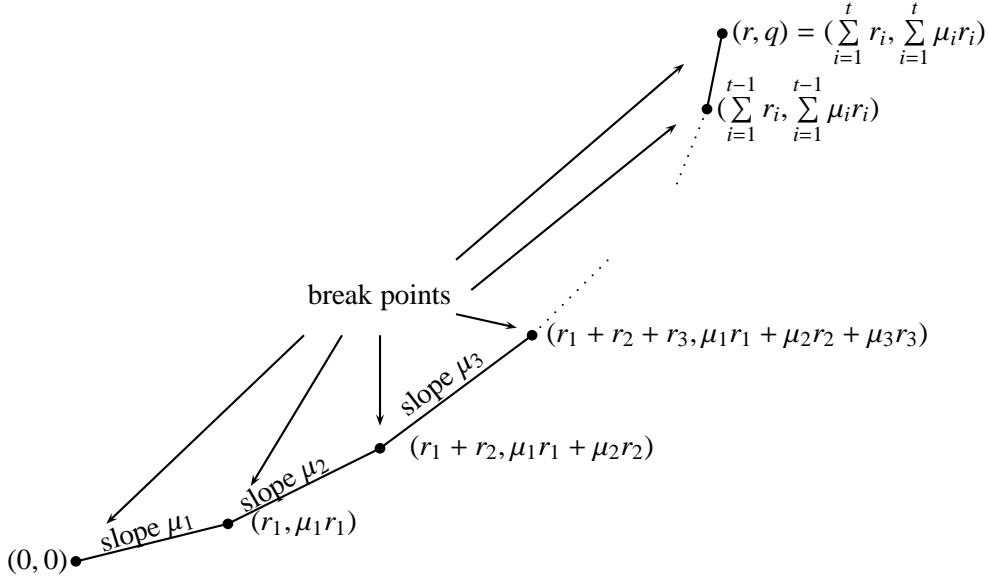
\begin{pspicture}(10,7)
\psdots[dotsize=4pt](0,0)(2,0.5)(4,1.5)(6,3)(8.5,7)(8.3,6)
\psline(0,0)(2,0.5)
\psline(2,0.5)(4,1.5) \psline(4,1.5)(6,3)
\psline[linestyle=dotted](6,3)(7,4)
\psline[linestyle=dotted](7.9,5)(8.3,6)
\psline{->}(3,3)(0.4,0.5)
\psline{->}(3.5,3)(2.3,1)
\psline{->}(4,3)(4,1.8)
\psline{->}(5,3.28)(5.8,3.1)
\psline{->}(5,4)(8,6.6)
\psline{->}(5.2,3.6)(7.9,5.8)
\psline(8.3,6)(8.5,7)
\rput[tr]{17}(1.5,0.75){\text{slope $\mu_1$}}
\rput[tr]{30}(3.2,1.5){\text{slope $\mu_2$}}
\rput[tr]{38}(5.2,2.9){\text{slope $\mu_3$}} 
\rput(4,3.5){break points}
\rput(-0.5,0){$(0,0)$}
\rput(3,0.5){$(r_1,\mu_1r_1)$}
\rput(6,1.5){$(r_1+r_2,\mu_1r_1+\mu_2r_2)$}
\rput(8.7,3){$(r_1+r_2+r_3,\mu_1r_1+\mu_2r_2+\mu_3r_3)$}
\rput(10.3,7){$(r,q)=(\sum\limits_{i=1}^t r_i,\sum\limits_{i=1}^t\mu_ir_i)$}
\rput(9.5,6){$(\sum\limits_{i=1}^{t-1} r_i,\sum\limits_{i=1}^{t-1}\mu_ir_i)$}
\end{pspicture} 

\vskip 1cm

Remark: From the definition of the Newton slopes, we have that all break points have integer coordinates.\\

Let $m\in \N$. Let $R$ be an $\F_p$-algebra, let $W_m(R)$ be the ring of Witt vectors of length $m$ with coefficients in $R$ and let $\sigma$ be the Frobenius endomorpisms of $W_m(R)$ for any $m \in \N$. Let $\M^n(W_m(R))$ be the abelian category whose objects are $W_m(R)$-modules endowed with $\sigma^n$-linear endomorphisms and whose morphisms are $W_m(R)$-linear maps that respect the $\sigma^n$-linear endomorphisms. We identify $\M^n(W_m(R))$ with a full subcategory of $\M^n(W_{m+1}(R))$. With a little abuse of terminology, we call the following morphism modulo $p^m$: $$\M^n(W_{m+s}(R))\to \M^n(W_m(R))$$ $$O\mapsto O\otimes_{W_{m+s}(R)}W_m(R),$$ for any $s\in \N_{>0}$. If $S$ is an $\F_p$-scheme, in a similar way we define $\M^n(W_m(S))$. If $S=\Spec~R$, we identify $\M^n(W_m(R))=\M^n(W_m(S))$.\\

In general, we can define an $F^n$-crystal $\C$ over any $\F_p$-algebra $R$, cf. [Ka79, (2.1)]. The evaluation of the $F^n$-crystal $\C$ at the thickening $(R\hookrightarrow W_m(R))$ is a triple $(_m\C,F,{_m\nabla})$, where $_m\C$ is a locally free $W_m(R)$-module of finite rank, $F:{_m\C}\to{_m\C}$ is a $\sigma^n$-linear endomorphism, and $_m\nabla$ is an integrable and topologically nilpotent connection on $_m\C$ that satisfies certain axioms. In this paper, connections as $_m\nabla$ will play no role. A morphism $\phi:\C\to\C_1$ of $F^n$-crystals over $R$ defines naturally a morphism in the category of $\M^n(W_m(R))$ $$_m\phi: {_m\C}\to {_m\C_1}.$$ The association $\phi\to{_m\phi}$ defines a $\Z_p$-linear (evaluation) functor from the category of $F^n$-crystals over $R$ into the category $\M^n(W_m(R))$.\\ 

Let $\O_1$ and $\O_2$ be two objects of $\M^n(W_m(R))$ such that their underlying $W_m(R)$-modules are locally free of finite ranks. Let $S=\Spec~R$. We consider the functor $$\textbf{Hom}(\O_1,\O_2):\mbox{Sch}^{S}\to\mbox{SET}$$ from the category Sch$^S$ of $S$-schemes to the category SET of sets, with the property that {\bf Hom}$(\O_1,\O_2)(S_1)$ is the set underlying the $\Z /p^m\Z$-module of morphisms of $\M(W_m(S_1))$ that are between $f^*(\O_1)$ and $f^*(\O_2)$; Here $f:S_1\to S$ is the structual morphism of $S_1$ and $f^{*}$ is the pullback to $S_1$

\begin{center}
$\xymatrix{f^{*}(\O_1)\ar[r]\ar[d]&\O_1\ar[d]\\
S_1\ar[r]& S}
\hskip 1in
\xymatrix{f^{*}(\O_2)\ar[r]\ar[d]&\O_2\ar[d]\\
S_1\ar[r]& S.}$
\end{center}

\textbf{Lemma 0} {\it The functor {\bf Hom}$(\O_1,\O_2)$ is representable by an affine $S$-scheme which locally is of finite presentation.}\\

{\bf Proof of Lemma 0:} In [Va06, (Lemma 2.8.4.1)], Vasiu proved the lemma in the case when $n=1$. We note that the proof for a general natural number $n$ goes the same way. $\blacksquare$\\ 

Now let $R$ be a reduced $\varmathbb{F}_p$-algebra. Let $R^{\perf}$ be its perfect closure.

Example: If $R=k[x_1,\cdots,x_n]$, then $R^{\perf}=\bigcup \limits_{u\ge 1} k[x_1^{\frac{1}{p^u}},\cdots,x_n^{\frac{1}{p^u}}]$.\\

Talking about an $F^n$-crystal over $S=\text{Spec}~ R$, we can look at its pullback to its perfect closure $R^{\perf}$. The pullback of such an $F$-crystal to $R^{\perf}$ is a finite rank locally free $W(R^{\perf})$-module $M$ equipped with a Frobenius linear map $F :M\to M$ such that $F (M)\supset p^t M$ for some $t\ge 0$.\\

Let $S=$Spec $R$ be a reduced affine $\varmathbb{F}_p$-scheme. Let $\mathfrak{C}$ be an $F^n$-crystal over $S$. For any $\varmathbb{F}_p$-homomorphism $\phi:R\to k$, the pullback $\mathfrak{C}^{(\phi)}$ is an $F^n$-crystal over $k$. Its Newton polygon depends only on the underlying point Ker$(\phi) \in S$. This allows us to speak about the Newton slopes and Newton polygons of $\mathfrak{C}$ at various points of $S$.\\

We define a function 
\begin{center}
$f_{\mathfrak{C}}:S\to$ Set of Newton Polygons\\
$s \mapsto$ NP($\mathfrak {C}_s$), where $\mathfrak {C}_s$ is the pullback of $\mathfrak{C}$ to the algebraic closure $\overline{k(s)}$ of the residue field $k(s)$ of the point $s\in S$, i.e., \\
$\xymatrix{
   \mathfrak{C}_s  \ar@{->}[r] \ar@{->}[d] & {\mathfrak{C}} \ar@{->}[d] \\
 {\mbox{Spec}~ \overline{k(s)}}\ar@{->}[r] & S.}$\\
  \end{center}

By the Newton polygon stratification of an $\varmathbb{F}_p$-scheme $S$ defined by an $F^n$-crystal $\mathfrak{C}$ over $S$ we mean the stratification of $S$ with the property that each stratum of it is of the form $f_{\mathfrak{C}}^{-1}(\mbox{a fixed Newton Polygon})_{\mbox{red}}$. For a given Newton polygon $\nu$,  let $S_{\nu}=f_{\mathfrak{C}}^{-1}(\nu)_{\mbox{red}}$.\\

Clearly we have a set partition: $S=\displaystyle\bigsqcup_{\nu\in I}S_{\nu}$, where $I$ is the set of all possible Newton Polygons.\\

Example: Let $R=k[t]$ and $R^{\perf}=k[t]^{\perf}$. Let $S=$ Spec $R$ be as above.  Let $\underline t\in W(R^{\perf})$ be the image of $(t,0,0,\cdots)\in W(R)$. Suppose $M$ is an $F$-crystal over $R$ such that its pullback to $R^{\perf}$ is an $F$-crystal $(M^{\perf}, F)$ where $M^{\perf}=W(R^{\perf})\oplus W(R^{\perf})$ with $\{e_1,e_2\}$ as a basis and $F(e_1)=\underline t e_1+pe_2$, $F(e_2)=pe_1$. We observe that when $t=0$ the corresponding Newton polygon $\nu_1$ has no break point in the middle (i.e., it is not the starting or ending point) and a unique slope $1$. When $t\neq 0$, the corresponding Newton polygon $\nu_2$ has the break point $(1,0)$ in the middle and two Newton slopes $0$ and $2$. Thus in this case as sets we have $S_{\nu_1}$ = $\{(t)\}$, $S_{\nu_2}=S\backslash \{(t)\}$ and $S_\nu=\emptyset$ for $\nu\neq \nu_1$ or $\nu_2$. \\


\textbf{Theorem 1 (Grothendieck--Katz, [Ka79, Theorem 2.3.1])} {\it Let $S$ be an $\varmathbb{F}_p$-scheme, $\C$ be an $F$-crystal over $S$ and $\nu$ be a Newton polygon. Then the set $$S_{\ge \nu}=\{s\in S|NP(\C_s)\mbox{ lies above }\nu\}$$ is Zariski closed in $S$.}\\

\textbf{Corollary 1} {\it Let $S$ be an $\varmathbb{F}_p$-scheme, $\C$ be an $F$-crystal over $S$ and $\nu$ be a Newton polygon. All the strata $S_\nu$ of the Newton polygon stratification of $S$ defined by $\C$ are locally closed subschemes of $S$.}\\

\noindent {\bf Types of purity notions:}\\

Let $S$ be an $\varmathbb{F}_p$-scheme. Let $T$ be a reduced locally closed subscheme of $S$. Let $\overline T$ be the schematic closure of $T$, i.e., topologically $\overline T$ is the Zariski closure of $T$ in $S$ and endowed with the reduced ringed structure.\\

(\textbf{\`a la Nagata--Zariski}) Suppose $S$ is locally Noetherian. We say $T$ is \textbf{weakly pure }in $S$, if each non-empty irreducible component of $\overline T\backslash T$ has pure codimention $1$ in $\overline T$. \\

Suppose $S$ is locally Noetherian. We say $T$ is\textbf{ universally weakly pure} in $S$, if for every locally Noetherian scheme $S_1$ equipped with a morphism $S_1\to S$, the locally closed subscheme $(T_{\times_S}S_1)_{\mbox{red}}$ is weakly pure in $S_1$.\\

(\textbf{Vasiu}) We say $T$ is \textbf{pure} in $S$, if $T$ is an affine $S$-scheme.\\

(\textbf{folklore}) We say $T$ is \textbf{strongly pure} in $S$, if locally in the Zariski topology of $\overline T$ there exists a global function $f$ on $\overline T$ such that $T=\overline T_f$ is the largest open subscheme of $\overline T$ over which $f$ is invertible.\\

\textbf{Lemma 1} {\it Let $S$ be a reduced locally Noetherian $\F_p$-scheme. If $T\subset S$ is locally closed and $S$-affine, then $\overline T\backslash T$ is either $\emptyset$ or of pure codimension $1$ in $S$. In particular, purity implies weak purity.}\\

{\bf Proof of Lemma 1:} Since this is a local statement, we can assume both $S=\Spec~ R$ and $T=\Spec~ A$ are reduced affine schemes. By replacing $S$ by $\overline T$, we can assume $T$ is open dense in $S$. If $T=S$, the statement is proved. Otherwise, let $x$ be a generic point of an irreducible component of $S\backslash T$ and let $d=\dim\O_{S,x}$. We need to show $d=1$. Since $T$ is dense in $S$ and $S\backslash T\neq \emptyset$, we know $d \ge 1$. Consider the pullback $\hat T$ of $T$ in the following commutative diagram:
\begin{center}
$\xymatrix{\hat T_{\red}\ar@{^{(}->}[r]\ar[d]&\Spec~ \widehat{\O_{S,x}}_{\red}\ar[d]\\
 T\ar@{^{(}->}[r]&S=\Spec~ R.}$
\end{center}

By replacing $S$ by $\Spec~\widehat{\O_{S,x}}_{\red}$ and $T$ with $\hat T_{\red}$, we can assume $R$ is a local, reduced, complete noetherian $\F_p-$algebra and all the points in $S\backslash T$ are of codimension $\ge d$ in $S$. As $R$ is local, complete ring, it is also excellent (cf. [Hi80, (34.B)]). Thus the normalization $S^n$ of $S$ is a finite $S$-scheme. Consider the following commutative diagram:

\begin{center}
$\xymatrix{ T^n\ar@{^{(}->}[r]\ar[d]&S^n=\Spec ~R^n\ar[d]\\
 T\ar@{^{(}->}[r]&S=\Spec ~R.}$
\end{center}

Since the morphism $S^n\to S$ is both finite and surjective, for any preimage of $x$ in $S^n$, say $\tilde x\in S^n$, we have $\dim \O_{S^n,\tilde x}=\dim \O_{S,x}$. By replacing $x$ by $\tilde x$, $T$ by $T^n$ and $S$ by $S^n$, we can assume both $T=\Spec~ A$ and $S=\Spec~ R$ are affine, normal, reduced, Noetherian $\F_p$-schemes and all the points in $S\backslash T$ are of codimension $\ge d$. For $1\le i\le n$, let $T_i=\Spec~A_i\subset S_i=\Spec~R_i$ be the irreducible components of $T$ and $S$ respectively and we have $A_i$ and $R_i$ are all normal, integral, Noetherian $\F_p$-algebras. Now Suppose $d\ge 2$. At the level of rings, the codimension $1$ points are height $1$ prime ideals. Therefore for any $1\le i\le n$ we have, $$R_i\hookrightarrow A_i=\displaystyle\bigcap_{p\mbox{ a prime of height }1}A_{i,p}=\displaystyle\bigcap_{p\mbox{ a prime of height }1}R_{i,p}=R_i $$ (cf. [Hi80, (17.H)] for the first and third equalities). Therefore $$S=\Spec\prod\limits_{i=1}^nR_i=\Spec\prod\limits_{i=1}^nA_i=T,$$ which is a contradiction. Hence $d=1$ and Lemma 1 is proved. $\blacksquare$\\

\noindent We have the following obvious implications and identifications:\\

Strong purity$\Rightarrow$ purity $\Rightarrow$ univeral weak purity $\Rightarrow$ weak purity

Purity=universal purity

Strong purity= universal strong purity\\

$\textbf{Theorem 2 (A. J. de Jong and F. Oort, [JO00, (Theorem 4.1)])}$ {\it Let $S$ be a reduced locally Noetherian $\varmathbb{F}_p$-scheme and let $\mathfrak{C}$ be an $F-$crystal over $S$. Then the Newton polygon stratification of $S$ defined by $\mathfrak{C}$ is universally weakly pure in $S$.} \\

$\textbf{Theorem 3 (A. Vasiu, [Va06, (Theorem 6.1)])}$ {\it Let $\mathfrak{C}$ be an $F$-crystal over a reduced locally Noetherian $\varmathbb{F}_p$-scheme $S$. Then the Newton polygon stratification of $S$ defined by $\mathfrak{C}$ is pure in $S$.}\\

Let $P_0$ be a point in the $xy$-coordinate plane and let $$S_{P_0}=\{s\in S|~NP(\mathfrak C_s)\mbox{ has }P_0\mbox{ as a break point}\}.$$ It can be shown that topologically $S_{P_0}$ is locally closed in $S$ (We will prove this in the proof of Theorem 5), and we endow it with the reduced ringed structure.\\

$\textbf{Theorem 4 (Y. Yang, [Ya10, (Theorem 1.1)])}$ {\it Let $S$ be a reduced locally Noetherian $\varmathbb{F}_p$-scheme and let $\mathfrak{C}$ be an $F$-crystal over $S$. Fix a point $P_0$ in the $xy$-coordinate plane. Then $S_{P_0}$ is universally weakly pure in $S$.}\\

\noindent Our main result is the following theorem, which will imply Theorems 2 to 4:\\

\textbf{Theorem 5 (J. Li)}  {\it Let $S$ be a reduced locally Noetherian $\varmathbb{F}_p$-scheme. Let $\mathfrak{C}$ be an $F^n$-crystal over $S$, $n\ge 1$. Fix a point $P_0$ in the $xy$-coordinate plane. Then $S_{P_0}$ is pure in $S$.}\\

{\bf Proof of Theorem 5:} It will be done in 5 steps. \\

\noindent\textbf{Step 1. Reduction step.}\\

Since purity is a local statement, we first assume $S=\Spec~ R$ is affine and we need to show that $S_{P_0}$ is an affine scheme. Let $P_0=(a,b)$. If $(a,b)\notin \N^2$ or $a=0$ and $b\neq 0$, then $S_{P_0}=\emptyset$ and the theorem holds trivially. If $(a,b)=(0,0)$, then $S_{p_0}=S$. Again the theorem holds trivially. Now we suppose $a,b\in \N_{>0}$. By replacing $\C$ by $\wedge^a \C$, we see that $S_{(a,b)}=\{s\in S|~NP(\mathfrak \wedge^aC_s)\mbox{ has } (1,b) \mbox{ as a break point}\}$. Therefore, we can assume $a=1$. By replacing $\C$ by $\C^{\otimes c}$ with $c$ a large integer (For example $c=r!$, where $r=$ rank($\C$)), we can assume that for each point $s\in S$, all Newton slopes of $\C_s$ are integers. Now let $S_{\ge\nu_1}=\{s\in S|NP(\C_s)\ge\nu_1\}$, where $\nu_1$ is the following Newton polygon:

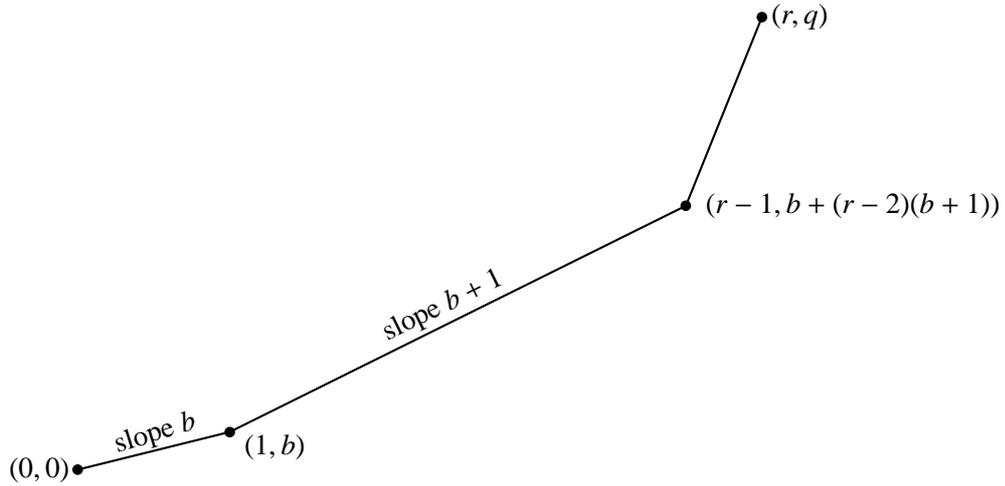
\begin{center}
\begin{pspicture}(10,7)
\psdots[dotsize=4pt](0,1)(2,1.5)(8,4.5)(9,7)
 \psline(0,1)(2,1.5)
\psline(2,1.5)(8,4.5) \psline(8,4.5)(9,7)

\rput[tr]{17}(1.5,1.8){\text{slope $b$}}
\rput[tr]{26}(5.5,3.7){\text{slope $b+1$}}

\rput(-0.5,1){$(0,0)$}
\rput(2.6,1.3){$(1,b)$}
\rput(10.2,4.5){$(r-1,b+(r-2)(b+1))$}
\rput(9.5,7){$(r,q)$}

\end{pspicture} 
\end{center}

From Theorem 1, $S_{\ge\nu_1}$ is closed in $S$, thus affine and since $(1,b)$ is a break point of $\nu_1$, we have $S_{P_0}\subset S_{\ge\nu_1}$. By replacing $S$ by $S_{\ge\nu_1}$, we can assume $NP(C_s)\ge \nu_1$ for any $s\in S$. Let $S_{\ge\nu_2}=\{s\in S|NP(\C_s)\ge\nu_2\}$, where $\nu_2$ is the following Newton polygon:

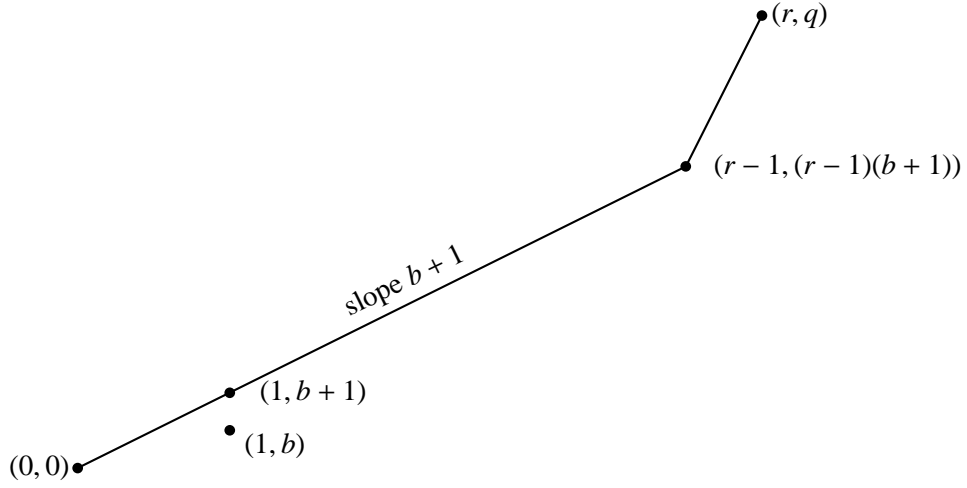
\begin{center}
\begin{pspicture}(10,7)
\psdots[dotsize=4pt](0,1)(2,1.5)(2,2)(8,5)(9,7)
 \psline(0,1)(8,5)
\psline(8,5)(9,7)
\rput[tr]{26}(5,4){\text{slope $b+1$}}
\rput(-0.5,1){$(0,0)$}
\rput(3.1,2){$(1,b+1)$}
\rput(2.6,1.3){$(1,b)$}
\rput(10,5){$(r-1,(r-1)(b+1))$}
\rput(9.5,7){$(r,q)$}
\end{pspicture} 
\end{center}

We have $S_{P_0}=S\backslash S_{\ge\nu_2}$ and this shows that $S_{P_0}$ is locally closed in $S$. If $ S_{\ge\nu_2}=\emptyset$, then $S_{P_0}=S$ is affine and the theorem is proved. Now we suppose $S_{\ge\nu_2}\neq \emptyset$ and $S_{P_0}=S\backslash S_{\ge\nu_2}$ is an open subscheme of $S=\Spec~ R$ and we need to show $S_{P_0}$ is affine. The statement is local in the faithfully flat topology of $S$ and thus we can assume that $S$ is local. Let $\hat R$ be the completion of $R$ and let $\hat S:=\Spec~\hat R$. As $\hat S$ is a faithfully flat $S$-scheme, to show that $S_{P_0}$ is affine it suffices to show that $S_{P_0}\times_{S}\hat S$ is affine. Let $\hat S_1=\Spec ~\hat R_1,\cdots,\hat S_j=\Spec ~\hat R_j$ be the irreducible components of the reduced scheme of $\hat S$ (Here $j\in \N$); They are spectra of local, complete, integral, Noetherian $\F_p$-algebras. The scheme $S_{P_0}\times_{S}\hat S$ is affine if and only if the irreducible components $S_{P_0}\times_{S}\hat S_1$, $\cdots$, $S_{P_0}\times_{S}\hat S_j$ of the reduced scheme of $S_{P_0}\times_{S}\hat S$ are all affine, cf. Chevalley's theorem in [Gr61, Ch. II, Cor. (6.7.3)]. So to prove the theorem we can assume $R=\hat R=\hat R_1$. As $R$ is a local, complete ring, it is also excellent, cf. [Hi80, (34.B)]. Thus the normalization $S^n$ of $S$ is a finite $S$-scheme. So $S^n$ is a semilocal, complete, integral, normal scheme. This implies $S^n$ is local. But $S_{P_0}$ is affine if and only if $S_{P_0}\times_{S} S^n$ is affine, cf. [Va06, (Lemma 2.9.2)]. Thus to prove the theorem, we can also assume $S$ is normal, i.e., $S=S^n$. We emphasize that for the rest of the proof we will use the fact that $R$ is a complete, integral, local, normal $\F_p$-algebra and $U:=S_{P_0}$ is an open subscheme of $S=\Spec ~R$. We can assume $U\neq \emptyset$ and hence it is dense in $S$. Let $k_S$ be the field of fractions of $R$ and $k$ be its algebraic closure.\\ 

\noindent \textbf{Step 2. The affine $\mathbf{S-}$scheme $\mathbf{H_m}$.}\\

Now let $\C_0$ be an $F^n-$crystal of rank $1$ and slope $b$ over $\F_p$, that is $\C_0=(\Z_p, p^b\sigma^n)$. Let $\C_{0,S}$ be its pullback to $S$:
\begin{center}
$\xymatrix{\C_{0,S}\ar[r]\ar[d]&\C_0\ar[d]\\
S\ar[r]&\Spec~ \F_p.}$
\end{center}
Let $m\in \N$ and $m>>0$. Let $_m\C_{0,S}$ be the evaluation of $\C_{0,S}$ at $W_m(S)=W_m(R)$ and $_m\C$ be the evaluation of $\C$ at $W_m(S)$. We view these evaluations as $W_m(R)$-modules equipped with $\sigma^n$-linear endomorphisms. By Lemma 0, the functor {\bf HOM}($_m\C_{0,S},_m\C$) is representable by an affine $S$-scheme $H_m$ which is of finite presentation and it is Noetherian since $S$ is Noetherian. In other words, we have an affine morphism of finite type: $H_m\to S$. Let $x$ be a point of codimension $1$ in $ U\xrightarrow{\mbox{open}}S$. Consider the stalk at $x$: $V_x:=\O_{U,x}=\O_{S,x}$. It is an integrally closed noetherian local domain of dimension one. Therefore, $V_x$ is a discrete valuation ring. Let $W_x$ be a discrete valuation ring with an algebraically closed residue field such that $V_x\hookrightarrow W_x$ and to simplify our notation, let $W:=W_x^{\perf}$ be its perfect closure. Consider the pullback of $\C$ to $W$:
 \begin{center}
$\xymatrix{\C_W \ar[r]\ar[d]&\C\ar[d]\\
 \Spec ~W\ar[r]&S.}$
\end{center}
As $S=S_{\mu_1}$, all Newton slopes of $\C_W$ are greater or equal to $b$ at every point of $\Spec~W$. By [Ka79, (Theorem 2.6.1)], we have an isogeny $\phi$: $\C_W \to \C'_W$, where $\C'_W$ is an $F^n-$crystal over $W$ which is divisible by $b$ and the cokernel of $\phi$ is annihilated by $p^t$ for $t=(r-1)b$. Let  $\varphi$ be the isogeny from $\C'_W$ to $\C_W$. Next we prove the following lemma:\\

\textbf{Lemma 2} {\it Let $W$ be a perfect discrete valuation ring of characteristic $p$ with an algebraically closed residue field. Let $(M,F)$ be an $F^n$-crystal over $W$ such that at each point of $\Spec~W$, all its Newton slopes are greater than or equal to $b$ and it is divisible by $b$. Assume the multiplicity of its Newton slope $b$ is $1$ at each point of $\Spec~W$, then there exists a unique sub-$F^n$-crystal $(M_1,F_1)$ of $(M,F)$ which is also a direct summand, where $(M_1,F_1)$ is an $F^n$-crystal over $W$ of rank $1$, Newton slope $b$.}\\

{\bf Proof of Lemma 2}: Since $(M,F)$ is divisible by $b$, we have $F(M)\subset p^bM$. Let $w$ be an arbitrary point in $\Spec~ W$ and let $(M_w,F_w)$ be the pullback of $(M,F)$ to $\overline{k(w)}$, where $\overline{k(w)}$ is the algebraic closure of the residue field $k(w)$ at $w$
$$\xymatrix{(M_w,F_w) \ar[r]\ar[d]&(M,F)\ar[d]\\
 \Spec ~\overline{k(w)}\ar[r]&\Spec~W.}$$
We have $F_w(M_w)\subset p^bM_w$ and $(1,b)$ is a break point of the Newton polygon of $(M_w,F_w)$ for every $w\in \Spec~W$. If $F_w(M_w)\subset p^{b+1}M_w$, then all the Newton slopes of $(M_w,F_w)$ will be greater than $b$, which is a contradiction. Thus $F_w(M_w)\nsubset p^{b+1}M_w$ and $F_w(M_w)\subset p^bM_w$. Since for every point $w\in \Spec~W$ the multiplicity of the Newton slope $b$ of $(M_w,F_w)$ is $1$, therefore the point $(1,b)$ lies on the Hodge polygon of $(M_w,F_w)$ at every point $w\in \Spec~W$ and by [Ka79, (Theorem 2.4.2)], lemma holds. $\blacksquare$\\

Applying Lemma 2, we have a monomorphism: $\tilde\C_{0,W}\hookrightarrow\C'_W$, where $\tilde\C_{0,W}$ is an $F^n-$crystal over $W$ of rank $1$ and Newton slope $b$ and $j_m$ admits a unique splitting. Modulo $p^m$, at the level of evaluation we get a monomorphism: $_m\tilde\C_{0,W}\hookrightarrow{_m\C'_W}$. Let $\C_{0,W}$ be the pullback of $\C_{0,S}$ to $\Spec ~W$:
\begin{center}
$\xymatrix{\C_{0,W} \ar[r]\ar[d]&\C_{0,S}\ar[d]\\
 \Spec~ W\ar[r]&S.}$
\end{center}

Since $\C_{0,W}$ is also an $F^n$-crystal over $W$ of rank $1$ and Newton slope $b$, modulo $p^m$ we have $_m\C_{0,W}$ isomorphic to $_m\tilde\C_{0,W}$. Now we have a morphism $i_W(m)$ from $_m\C_{0,W}\to$ $_m\C_W$ by composing the following morphisms: $_m\C_{0,W}\simeq$ $_m\tilde \C_{0,W}\hookrightarrow$ $_m\C'_W\xrightarrow{\varphi_m}$$_m\C_W$, where $\varphi_m$ is the isogeny $\varphi$ modulo $p^m$, and thus its cokernel is annihilated by $p^t$. Let $j_m$ be the composition of the morphisms:  $_m\C_{0,W}\simeq$ $_m\tilde \C_{0,W}\hookrightarrow$ $_m\C'_W$ and we know that $j_m$ is a monomorphism that splits.\\

\noindent \textbf{Step 3. Gluing morphisms.}\\
 
Before we glue the morphisms, let us first discuss three useful cases of inductive limits.\\

Let $V\hookrightarrow V_1$ be a monomorphism of commutative $\F_p$-algebras. Suppose we have an inductive limit $V_1$ = ind $\lim\limits_{\alpha\in \Lambda} V_{\alpha}$ of commutative $V$-subalgebras of $V_1$ indexed by the set of objects $\Lambda$ of a filtered, small category. For $\alpha\in \Lambda$, let $f^{\alpha} : \Spec~V_\alpha\to \Spec~V$ be the natural morphism. Let $(O,\phi_O)$ and $(O',\phi_{O'} )$ be objects of $\M^n(W_m(V ))$ such that $O$ and $O'$ are free $W_m(V )$-modules of finite rank. Let $(O_1,\phi_{O_1})$ and $(O'_ 1,\phi_{O'_1})$ be the pullbacks of $(O,\phi_O)$ and $(O',\phi_{O'})$ (respectively) to objects of $\M^n(W_m(V_1))$. We consider a morphism $$u_1 : (O_1,\phi_{O_1} )\to(O'_1,\phi_{O'_1})$$ of $\M^n(W_m(V_1))$. We fix ordered $W_m(V)$-bases $\B_O$ and $\B_O'$ of $O$ and $O'$ (respectively). Let $B_1$ be the matrix representation of $u_1$ with respect to the ordered $W_m(V_1)$-basis of $O_1$ and $O'_1$ defined naturally by $\B_O$ and $\B_O'$ (respectively). Let $V_{u_1}$ be the $V$-subalgebra of $V_1$ generated by the components of the Witt vectors of length $m$ with coefficients in $V_1$ that are entries of $B_1$.  As $V_{u_1}$ is a finitely generated $V$-algebra, there exists $\alpha_0\in \Lambda$ such that $V_{u_1}\hookrightarrow V_{\alpha_0}$. This implies that $u_1$ is the pullback of a morphism $$u_{\alpha_0} : f^{\alpha_0\ast}_m (O,\phi_O)\to f^{\alpha_0\ast}_m (O' ,\phi_{O'})$$ of $\M^n(W_m(V_{\alpha_0})).$ Here are three special cases of interest.\\
\begin{itemize}
\item[(a)] If $V$ is a field and $V_1$ is an algebraic closure of $V$ , then as $V_\alpha\mbox{'}$s we can take the finite field extensions of $V$ that are contained in $V_1$.
\item[(b)] If $V_1$ is a local ring of an integral domain $V$ , then as $V_\alpha\mbox{'}$s we can take the $V$-algebras of global functions of open, affine subschemes of $\Spec~V$ that contain $\Spec~V_1$.
\item[(c)] We consider the case when $V$ is a discrete valuation ring that is an $N-2$ ring in the sense of [Hi80, (31.A)], when $V_1$ is a faithfully flat $V$-algebra that is also a discrete valuation ring, and when each $V_\alpha$ is a $V$-algebra of finite type. The flat morphism $f^{\alpha_0} : \Spec~V_{\alpha_0}\to \Spec~V$ has quasi-sections, cf. [Gr64, Ch. IV, Cor. (17.16.2)]. In other words, there exists a finite field extension $\tilde k$ of $k$ and a $V$-subalgebra $\tilde V$ of $\tilde k$ such that: (i) $\tilde V$ is a local, faithfully flat $V$-algebra of finite type whose field of fractions is $\tilde k$, and (ii) we have a morphism $\tilde f^{\alpha_0} : \Spec~\tilde V\to \Spec~V_{\alpha_0}$ such that $\tilde f := f^{\alpha_0}\cdot \tilde f^{\alpha_0}$ is the natural morphism $\Spec~\tilde V\to \Spec~V$. As $V$ is an $N-2$ ring, its normalization in $\tilde k$ is a finite $V$-algebra and so a Dedekind domain. This implies that we can assume $\tilde V$ is a discrete valuation ring. For future use, we recall that any excellent ring is a Nagata ring (cf. [Hi80, (34.A)]) and so also an $N-2$ ring (cf. [Hi80, (31.A)]). Let $$\tilde u:\tilde f^\ast_m(O,\phi_O)=\tilde f^{\alpha_0\ast}_m(f^{\alpha_0\ast}_m(O,\phi_O))\to\tilde f^{\ast}_m(O',\phi_{O'})=\tilde f^{\alpha_0\ast}_m(f_m^{\alpha_0\ast}(O',\phi_{O'}))$$
be the pullback of $u_{\alpha_0}$ to a morphism of $\M^n(W_m(\tilde V))$. If $V$ is the local ring of an integral $\F_p$-scheme $U$, then $\tilde V$ is a local ring of the normalization of $U$ in $\tilde k$. So from (b) we get that there exists an open subscheme $\tilde U$ of this last normalization that has $\tilde V$ as a local ring and that has the property that $\tilde u$ extends to a morphism of $\M^n(W_m(\tilde U))$.\\
\end{itemize}
\newpage
\noindent Facts:
\begin{itemize}
\item[1.] If $u_1$ is a monomorphism and $(O_1,\phi_{O_1})$ is a direct summand of $(O'_ 1,\phi_{O'_1})$, then $u_{\alpha_0}$ is a monomorphism and  $ f^{\alpha_0\ast}_m (O,\phi_O)$ is a direct summand of $f^{\alpha_0\ast}_m (O' ,\phi_{O'})$.
\item[2.] If $u_1$ is a morphism such that its cokernel is annihilated by $p^t$, by enlarging $V_{u_1}$, we can assume Coker($u_{\alpha_0}$) is also annihilated by $p^t$, cf. [Va06, 2.8.3].
\end{itemize}

Now let $$v:=\mbox{ max }\{v(1,1,b,c)~|~c=0,1,2,\cdots,\mbox{ Maximum hodge slope of }\C_k\},$$ cf. [Va06, 5.1.1(b)] for the function $v(\cdot,\cdot,\cdot,\cdot)$ with $M_1=\C_{0,k}$ and $M_2=\C_k$ ($v$ does not depend on $m$). Replacing $m$ by $m+v$, from the above discussion (case (c)), we get that there exists a finite field extension $k_{S,\tilde V_x}$ of $k_S$ and an open, affine subscheme $U_{\tilde V_x}$ of the normalization of $U$ in $k_{S,\tilde V_x}$, such that $U_{\tilde V_x}$ has a local ring $\tilde V_x$ which is a discrete valuation ring that dominates $V_x$ and moreover we have a morphism
\begin{center}
$i_{U_{\tilde V_x}}(m+v):{_{m+v}\C}_{0,U_{\tilde V_x}}\to{_{m+v}\C}_{U_{\tilde V_x}},$
\end{center}
 where $_{m+v}\C_{0,U_{\tilde V_x}}$ and $_{m+v}\C_{U_{\tilde V_x}}$ are the pullbacks of $_{m+v}\C_{0,S}$ and $_{m+v}\C$ to $\Spec~ U_{\tilde V_x}$ respectively. Modulo $p^m$, we have a morphism \begin{center}
$i_{U_{\tilde V_x}}(m):{_{m}\C}_{0,U_{\tilde V_x}}\to{_{m}\C}_{U_{\tilde V_x}},$
\end{center} where $_{m}\C_{0,U_{\tilde V_x}}$ and $_{m}\C_{U_{\tilde V_x}}$ are the pullbacks of $_{m}\C_{0,S}$ and $_{m}\C$ to $\Spec~ U_{\tilde V_x}$ respectively. Let $I_m$ be the set of morphisms $_m\C_{0,k}\to$ $_m{\C_k}$ that are reductions modulo $p^m$ of morphisms $_{m+v}\C_{0,k}\to$ $_{m+v}{\C_k}$. From [Va06, (Theorem 5.1.1(b)+ Remark 5.1.2)], we get that each morphism in $I_m$ lifts to a morphism $\C_{0,k}\to\C_k$; Thus $I_m$ is a finite set. Based on the above discussion (case (a)), by replacing $(S,U)$ by its normalizations $(\tilde S,\tilde U)$ in a finite field extension of $k_S$, we can assume that $I_m$ is the set of pullbacks of a set of morphisms $L_m$ of $\M^n(W_m(k_S))$. Since $i_{U_{\tilde V_x}}(m)$ is the reduction modulo $p^m$ of the morphism $i_{U_{\tilde V_x}}(m+v)$, thus $i_{U_{\tilde V_x}}(m)\in I_m$ and the pullback of $i_{U_{\tilde V_x}}(m)$ to a morphism of $\M^n(W_m(k_{S,V_x}))$ is also the pullback of a morphism in $L_m$. As $V_x=\tilde V_x \cap k_S$, inside $W_m(k_{S,V_x})$ we have $W_m(V_x)=W_m(\tilde V_x)\cap W_m(k_S)$. This implies that the pullback of $i_{U_{\tilde V_x}}(m)$ to a morphism of $\M^n(W_m(\tilde V_x))$ is in fact the pullback of a morphism of $\M^n(W_m(V_x))$. From the above discussion (case (b)) (applied with $(V_1,V)$ replaced by $(V_x,R)$), we get the existence of an open subscheme $U_{V_x}$ of $U$ that has $V_x$ as a local ring and such that we have a morphism:
\begin{center}
$i_{U_{V_x}}(m):$ $_m\C_{0,U_{V_x}}\to$$_m\C_{U_{V_x}},$
\end{center}
where $_m\C_{0,U_{V_x}}$ and $_m\C_{U_{V_x}}$ are the pullbacks of $_m\C_{0,S}$ and $_m\C$ to $\Spec~ U_{V_x}$ respectively. Using the above Facts 1 and 2, if we apply similar arguments on $j_m:{_m\C}_{0,W}\hookrightarrow{_m\C}'_W$ and $\varphi_m:{_m\C}'_W\to{_m\C}_W$, we can assume that $i_{U_{V_x}}(m)$ is the composition of two morphisms $j_{U_{V_x}}(m):{_m\C}_{0,U_{V_x}}\to{_m\C'_{U_{V_x}}}$ and $\varphi_{U_{V_x}}(m):{_m\C'_{U_{V_x}}}\to {_m\C_{U_{V_x}}}$, where $j_{U_{V_x}}(m)$ is a monomorphism that splits and the cokernel of $\varphi_{U_{V_x}}(m)$ is annihilated by $p^t$.  Now at the level of modules consider the following diagram:
\begin {center}$\xymatrix{_m\C_{0,U_{V_x}}\ar[r]^{i_{U_{V_x}}(m)}\ar@{_{(}->}[d]& _m\C_{U_{V_x}}\ar@{_{(}->}[d]\\
_m\C_{0,k} \ar[r]^{i_{k,x}(m)}& _m\C_k.}$\end{center}

Let $_m\C_{0,k}=(W_m(k),p^b\sigma^n)$. Since the Newton polygon of $\C_k$ has $(1,b)$ as a break point, at the level of evaluation there is a unique rank $1$ free sub-$W_m(k)$-module $_m\tilde\C_{k}=(W_m(k),p^b\sigma^n)$ of $_m\C_k$. The morphism $i_{U_{V_x}}(m)$ is nonzero mod $p^t$ (assuming $m>t$) since $i_{U_{V_x}}(m)=\varphi_{U_{V_x}}(m)\circ j_{U_{V_x}}(m)$, where Coker($\varphi_{U_{V_x}}(m)$) is annihilated by $p^t$ and $j_{U_{V_x}}(m)$ is a monomorphism that splits. Therefore the morphism $i_{k,x}(m)$ is not a zero morphism if $m>t$, thus it must factor through $_m\tilde\C_{k}$ $$i_{k,x}(m):{_m\C}_{0,k}\xrightarrow{f_x(m)}{_m\tilde\C_{k}}\hookrightarrow{_m\C}_k.$$ Let $\epsilon=f_x(m)(1)\in W_m(k)$ and consider the following commutative diagram
\begin{center}
$\xymatrix{ _m\C_{0,k}=W_m(k) \ar[r]^{f_x(m)}\ar[d]_{p^b\sigma ^n}& _m\tilde\C_{k}=W_m(k)\ar[d]^{p^b \sigma ^n}\\
_m\C_{0,k}=W_m(k) \ar[r]^{f_x(m)} & _m\tilde \C_{k}=W_m(k).}$
\end{center}
We have $p^b\sigma^n(\epsilon)=p^b\sigma^n(f_x(m)(1))=f_x(m)(p^b\sigma^n(1))=f_x(m)(p^b)=p^bf_x(m)(1)=p^b\epsilon$. Therefore, $\sigma^n(\epsilon)=\epsilon$ mod $p^{m-b}$. Let $\epsilon=p^{t_x}\epsilon '$, where $\epsilon '\in W_m(k)$ is invertible. Since the cokernel of $f_x(m)$ is annihilated by $p^t$ (as the pullback to $\Spec~k$ of $\varphi_{U_{V_x}}(m)$ is annihilated by $p^t$ and the pullback to $\Spec~k$ of $j_{U_{V_x}}(m)$ splits), we have $1\le t_x\le t$. We now have $\sigma^n(\epsilon ')=\epsilon '$ mod $p^{m-b-t_x}$. Therefore $\sigma^n((\epsilon ')^{-1})=(\epsilon ')^{-1}$ mod $p^{m-b-t_x}$. This shows that multiplication by $(\epsilon ')^{-1}$ is an automorphism of $W_{m-b-t_x}(k)$ and of $W_{m-b-t}(k)$ since $t_x\le t$. Replacing $m$ by $m-b-t$ and composing $f_x(m)$ with the endomorphism of $W_m(k)$ $$ W_m(k)\xrightarrow{\mbox{multiplication by }(\epsilon ')^{-1}} W_m(k)$$ we can assume $f_x(m)(1)=1\in W_m(k)$.\\

Now let $y\in U$ be another point of codimension $1$. Similarly we can construct a morphism $f_y(m)$ and by the above construction $f_x(m)$ and $f_y(m)$ coincide. This tells us that the pullbacks of $i_{U_{V_x}}(m)$ and $i_{U_{V_y}}(m)$ to morphisms of $\M^n(W_m(U_{V_x}\cap U_{V_y}))$ coincide and therefore they glue together. Now we glue the morphisms $i_{U_{V_x}}(m)$ for all $x\in U$ with codimension $1$ and we obtain a morphism $i_{U_0}(m)$: $_m\C_{0,U_0}\to{_m\C}_{U_0}$, where $U_0$ is an open subscheme of $U$ and $S$ and codim $(U\backslash U_0)\ge 2$. By the definition of $H_m$, we have constructed an $S$-section $j$: $U_0\to H_m$.\\

\noindent{\bf Step 4. Sections of $\mathbf{H_m}$.}\\

By the gluing argument of Step 3, we have an open subscheme $U_0$ of $S$ with codim($U\backslash U_0)\ge 2$ and a section $j:~U_0\to H_m$ such that the following diagram commutes:

 \begin{center}
$\xymatrix{&&H_m\ar[d]\\
 U_0\ar@{^{(}->}[r]\ar[urr]^{\mbox{section }j}&U\ar@{^{(}->}[r]&S.}$
\end{center}

Let $J_m$ be the schematic closure $\overline{j(U_0)}$ of $j(U_0)$ in $H_m$. As $H_m$ is affine and Noetherian, $J_m$ is also an affine, Noetherian $S$-scheme. Since $S$ is Noetherian, normal and integral, we get that $U_0$ is also Noetherian, normal and integral and therefore $J_m$ is integral. Now we have the following commutative diagram:

 \begin{center}
$\xymatrix{&&J_m\ar[d]^{\mbox{affine}}\\
 U_0\ar@{^{(}->}[r]\ar[urr]^{\mbox{open}}&U\ar@{^{(}->}[r]&S.}$
\end{center}

Now consider the pullback $\tilde{J_m}$ of $J_m$ to $U$:

 \begin{center}
$\xymatrix{&\tilde{J_m}\ar@{^{(}->}[r]^{\mbox{open}}\ar[d]_{g}^{\mbox{affine}}&J_m\ar[d]^{\mbox{affine}}\\
 U_0\ar@{^{(}->}^{\mbox{open}}[r]\ar@{^{(}->}[ur]^{\mbox{open}}&U\ar@{^{(}->}[r]^{\mbox{open}}&S.}$
\end{center}

Claim: The morphism $g$ is an isomorphism, i.e., $\tilde J_m\simeq U$.\\

Proof of Claim: To prove that $g$ is an isomorphism, we can assume $U=\Spec ~A$ is affine. As $g$ is an affine morphism, we can also assume $\tilde J_m=\Spec~ B$ is also affine. Since $U_0$ is open dense in both $U$ and $\tilde J_m$, therefore $U$ and $\tilde J_m$ are birationally equivalent. Thus their fractional fields $\Frac(A)$ and $\Frac(B)$ are equal. As $\mbox{codim }(U\backslash U_0)\ge 2$ and the fact that $U_0$ is in both $U$ and $\tilde J_m$, we have $A_p=B_p$ for any prime $p\in \Spec~A$ of height $1$. As $A$ is a Noetherian normal domain, we have $$A\hookrightarrow B\subset\displaystyle\bigcap_{q\in\Spec ~B\mbox{ of height }1}B_q\subset \displaystyle\bigcap_{p\in \Spec~ A\mbox{ of height }1}A_p=A$$ (cf. [Hi80, (17.H)] for the equality part; The first monomorphism is given by $g^\#$). Therefore $A=B$ and the claim is proved.\\

Now we have a section from $U\simeq \tilde J_m\hookrightarrow H_m$.\\

\noindent{\bf Step 5. Final output:} $\mathbf{U=J_m \mbox{\textbf{ for }}m>>0.}$\\

{\bf Claim:} {\it If $m>>0,$ then $U=J_m$.}\\

{\bf Proof of Claim:} Suppose not. Let $\eta$ be the generic point of an irreducible component of $J_m\backslash U$. Consider the stalk at $\eta$: $\O_{J_m,\eta}:=R_p$, which is a normal, local Noetherian ring but not a field since $U$ is dense in $J_m$. Now we have $\dim(R_p)\ge 1$. Therefore, we can find a prime ideal $q\in U\subset S_p$, such that $\dim(R_p/q)=1$. Now the normalization ${(R_p/q)}^n$ of $R_p/q$ is an integrally closed noetherian local ring of dimension one, thus a discrete valuation ring. The morphism $\Spec~{(R_p/q)}\to J_m$ obtained from composing the natural morphisms $\Spec~{(R_p/q)}^n\to \Spec~R_p/q\to \Spec~R_p\to J_m$ sends the generic point of $\Spec~(R_p/q)^n$ to $U$ and the closed point to $J_m\backslash U$. With the help of a little commutative algebra, we can further assume the generic point of $\Spec~(R_p/q)^n$ is mapped into $U_0\subset U$. Let $D$ be the completion of the discrete valuation ring $(R_p/q)^n$, it is isomorphic to a power series ring $l[[t]]$ for some field $l$ of characteristic $p$ by Cohen structure theorem, cf. [Ei07, (Theorem 7.7)]. By injecting $l[[t]]$ to $\overline l[[t]]$, we can further assume $l$ is algebraically closed. Now we look at the pullback of $\C$ to $\Spec~ D$:
 \begin{center}
$\xymatrix{\C_D\ar[r]\ar[d]&\C\ar[d]\\
 \Spec~D\ar[r]&S.}
$\end{center}
It is an $F^n$-crystal such that at the generic point $\Spec~(\Frac(D))$ of $D$ its Newton polygon has $(1,b)$ as a break point and at the closed point $\Spec~l$ of $D$ all Newton slopes are at least $b+1$. Suppose the generic point $\Spec~(\Frac(D))$ is mapped to $z\in U_0$, we pull back the following morphisms to $\Spec~(D)$:
\begin{center}
$\xymatrix{ {_m\C}_{0,U_{V_z}}\ar[r]^{j_{U_{V_z}}(m)}\ar[d]&{_m\C'_{U_{V_z}}}\ar[d]\ar[r]^{\varphi_{U_{V_z}}(m)}&{_m\C}_{U_{V_z}}\ar[d]\\
{_m\C}_{0,D}\ar[r]&{_m\C'_D}\ar[r]&{_m\C}_D.}$ 
\end{center}

Let $E=D^{\perf}$ and $_m\C_{0,E}$, $_m\C_E$ and $_m\C'_E$ be the pullback of $_m\C_{0,D}$, $_m\C_D$ and $_m\C'_D$ to Spec $E$ respectively:
 \begin{center}
$\xymatrix{_m\C_{0,E}\ar[r]\ar[d]&_m\C_{0,D}\ar[d]\\
 \Spec~E\ar[r]&\Spec~D}$
\hskip 0.5in
$\xymatrix{_m\C_E\ar[r]\ar[d]&_m\C_D\ar[d]\\
 \Spec~E\ar[r]&\Spec~D}$
\hskip 0.5in
$\xymatrix{_m\C'_E\ar[r]\ar[d]&_m\C'_D\ar[d]\\
 \Spec~E\ar[r]&\Spec~D.}$
\end{center}
Recall $r=$ rank($\C$)=rank($_m\C$)= rank($_m\C_E$)= rank($_m\C'_E$). Let $_m\C_E=(W_m(E)^r,~\bar F)$, let $_m\C'_E=(W_m(E)^r,~\bar F')$ and let $_m\C'_{0,E}=(W_m(E),~p^b\sigma^n)$. Since $\C'$ is divisible by $b$, we can assume $\bar F'=p^b\bar G$. Let $x$ be a basis element of $M_0/p^mM_0$ (We can assume $\sigma^n x=x$.) and consider the following morphism at the level of $W_m(E)$-modules:
\begin {center}
${_m\C'}_{0,E}\xrightarrow{\gamma}{_m\C'_E}\xrightarrow{\gamma '} {_m\C}_E$\\
$x \mapsto  y\mapsto \tilde y$,
\end{center}
where $y=(y_1,y_2,\cdots,y_r)$ and $y_i\in W_m(E)$ for $i=1,2,\cdots,r$. We have $ p^b(\overline G(y)-y)=p^b\overline G (\gamma (x))-p^b\gamma(x)=\overline F'(\gamma (x))-\gamma(p^b x)=\gamma(p^b \sigma^n x)-\gamma(p^b x)=0$. Therefore $\overline G(y)=y$ mod $p^{m-b}$. Since the cokernel of $\gamma ': {_m\C'_E}\to {_m\C_E}$ is annihilated by $p^t$ and $\gamma:{_m\C'}_{0,E}\to{_m\C'_E}$ is a monomorphism that splits, we can write $y=p^{t_0}z$ for $0\le t_0\le t$ and $z=(z_1,z_2,\cdots,z_r)$ is not divisible by $p$, where $z_{i}=(z_{i,0},z_{i,1}z_{i,2},\cdots,z_{i,m-1})$ with $z_{i,j}\in E$  for $1\le i\le r,0\le j\le m-1$.\\

{\bf Subclaim 1:} {\it $y\neq 0$ mod $p^{t_0}$ at the closed point $\Spec~l$ of $\Spec~E$.}\\

{\bf Proof of Subclaim 1:} Recall $E=l[[T]]^{\perf}$. As $\overline G(y)=y$ mod $p^{m-b}$ and $y=p^{t_0}z$ with $0\le t_0\le t$, we have $\overline G(z)=z$ mod $p^{m-b-t}$. Since $z=(z_1,z_2,\cdots,z_r)$ is not divisible by $p$, modulo $p$ in $E^r$ we have $(z_{1,0},z_{2,0},\cdots,z_{r,0})\neq 0$. Suppose $z=0$ at the closed point $\Spec~l$, then for some $v\in \N$, we have $z_{i,0}\in l[[T^{\frac{1}{p^v}}]]\subset l[[T]]^{\perf}$ and $z_{i,0}=0$ mod $T^{\frac{1}{p^v}}$ for all $i=1,2,\cdots,r$. Now Let $\tilde T=T^{\frac{1}{p^v}}$ and let $z_{i,0}=\tilde T^f z'_{i,0}$ for some $f\in \N$, for all $i=1,2,\cdots, r$ and $z'_{i_0,0}\neq 0$ mod $\tilde T$ for some $1\le i_0\le r$. Let $\overline{\overline G}$ be $\overline G$ mod $p$. We have $\tilde T^f(z'_{1,0},z'_{2,0},\cdots,z'_{r,0})=(z_{1,0},z_{2,0},\cdots,z_{r,0})=\overline {\overline G}((z_{1,0},z_{2,0},\cdots,z_{r,0}))=\overline {\overline G}(\tilde T^f(z'_{1,0},z'_{2,0},\cdots,z'_{r,0}))=\tilde T^{p^nf}\overline {\overline G}((z'_{1,0},z'_{2,0},\cdots,z'_{r,0}))$. This contradicts to the fact that $z'_{i_0,0}\neq 0$ mod $\tilde T$ and thus $z\neq 0$ at the closed point. As $y=p^{t_0}z$, we have $y\neq 0$ mod $p^{t_0}$ at the closed point and Subclaim $1$ is proved.\\

Evaluating the morphism $\gamma$ at the closed point $\Spec~ l$ of $\Spec~E$ we have a morphism
\begin{center}$\beta:(W_m(l),p^b\sigma^n)\to (N,\hat F)$,\end{center}
where $(N,\hat F)=(\C_l$ mod $p^m$) and $\C_l$ is the pullback of $\C_D$ to the closed point $\Spec~l$. Since $y\neq 0$ mod $p^{t_0}$ at $\Spec~l$, the morphism $\beta$ is non-zero when reduced modulo $p^{t_0}$.\\

{\bf Subclaim $2$:} {\it $\C_l$ has Newton slope $b$.}\\

{\bf Proof of Subclaim 2:} Let $\C_l=(M_l,F_l)$ and suppose all Newton slopes of $\C_l$ are at least $b+1$. By [Ka79, (Sharp Slope Estimate 1.5.1)] we have $F_l^u(M_l)\subset p^{(b+1)u-h}M_l$ for all $u\in\N$, where $h>0$ is a fixed number. If $u>h+t_0+1$, we have $F_l^u(M_l)\subset p^{(b+1)u-h}M_l\subset p^{bu+t_0+1}M_l$. Let $(N,\hat F^u)=(M_l,F_l^u)\mbox{ mod }p^m$ and we also have $\hat F^u(N)\subset p^{bu+t_0+1}N$. Now consider the following morphism
\begin{center}$\alpha:(W_m(l),(p^b\sigma^n)^u)\to (N,\hat F^u)$\\
$1\mapsto e$,
\end{center}
where as a function $\alpha=\beta$. Since the morphism $\beta$ is non-zero when reduced modulo $p^{t_0}$, so is the morphism $\alpha$. However, 
$p^{bu}e=p^{bu}\alpha(1)=\alpha(p^{bu})=\alpha((p^b\sigma^n)^u(1))=\hat F^u(\alpha(1))=\hat F^u(e)\in p^{bu+t_0+1}N$ and thus $p^{m-t_0}\alpha(1)=p^{m-t_0}e=p^{m-t_0-bu}p^{bu}e\in p^{m-t_0-bu}p^{bu+t_o+1}N=p^{m+1}N=0$ if $m\ge bu+t+1\ge bu+t_0+1$, which is a contradiction since $\alpha$ mod $p^{t_0}$ is non-zero. Subclaim $2$ is therefore proved.\\

By Subclaim $2$, the $F^n$-crystal $\C_D$ has Newton slope $b$ at the closed point $\Spec~l$, a contradiction by our assumption. Therefore the claim is proved, i.e., $U=J_m$.\\

Now we have $S_{P_0}=U=J_m$ and $J_m$ is affine, therefore $S_{P_0}$ is affine and Theorem 5 is thus proved. $\blacksquare$\\

\noindent First application of Theorem 5:\\

\textbf{Proposition 1} {\it Theorem 5 implies Theorem 3, i.e., if $S_{p_0}$ is pure in $S$ for each point $P_0$ in the $xy$-coordinate plane, then the Newton polygon stratification of $S$ defined by $\mathfrak{C}$ is pure in $S$.}\\

{\bf Proof of Proposition 1:} Let $S$ and $\C$ be as above and let $\nu_0$ be a Newton polygon. We need to show that the set $$S_{\nu_0}=\{s\in S| NP(\C_s)=\nu_0\}$$ is an affine $S$-scheme. To prove that $S_{\nu_0}$ in $S$ is affine, we can further assume $S$ is affine. Let the break points of $\nu_0$ be $Q_0, Q_1,\cdots, Q_t$ and let $$S_{\ge \nu_0}=\{s\in S|NP(\C _s)\ge \nu_0\}.$$ Notice for any $s\in S$, the Newton polygon $NP(\C_s)$ is $\nu_0$ if and only if $NP(\C_s)\ge \nu_0$ and moreover $Q_0, Q_1,\cdots, Q_t$ are all break points of $NP(\C_s)$. Therefore $S_{\nu_0}=S_{\ge \nu_0}\cap S_{Q_0}\cap S_{Q_1}\cap\cdots\cap S_{Q_t}$. By Theorem 1, $S_{\ge \nu_0}$ is closed in $S$, thus affine. By Theorem 5, $S_{Q_i}$ is affine for $i=0,1,\cdots,t$. Since $S$ is affine thus separated, we conclude that $S_{\nu_o}$ is affine. $\blacksquare$\\

\noindent Second application of Theorem 5:\\

\textbf{Proposition 2} {\it Theorem 5 implies Theorem 4.}\\

{\bf Proof of Proposition 2:} From Theorem 5, we know that $S_{P_0}$ is pure. As purity implies universal weak purity, we have $S_{P_0}$ being universally weakly pure. $\blacksquare$\\

\noindent {\bf Now we have the following implications:}\\

Theorem 5 $\Rightarrow$ Theorem 4 $\Rightarrow$ Theorem 2. \\

Theorem 5 $\Rightarrow$ Theorem 3 $\Rightarrow$ Theorem 2.\\

\noindent Third application of Theorem 5:\\

Let $\C=(M,F)$ be an $F^n$-crystal over an algebraically closed field $k$ of char $p>0$. Let $\overline M$ be the reduction modulo $p$ of $M$ and let $\overline F:\overline M\to\overline M$ be the reduction modulo $p$ of $F$. Then the $p$-rank $t$ of $\C$ can be defined equivalently as follows:\\

(i) It is $\dim_{\F_p}(\{x\in\overline M|\overline F(x)=x\}).$\\

(ii) It is the multiplicity $t$ of the Newton slope $0$ of $\C$.\\

(iii) It is the unique non-negative integer such that $(t,0)$ is a break point of the Newton polygon of $\C$.\\

Let $S$ be a locally Noetherian $\F_p$-scheme and let $\C$ be an $F^n$-crystal over $S$. For each $t\in \N$, let $Y_t$ be the reduced, locally closed subscheme of $S$ formed by those points $s\in S$ with the property that the $p$-rank of the pullback of $\C$ to $\overline{k(s)}$, where $\overline{k(s)}$ is the algebraic closure of the residue field of $s$, is exactly $t$. We call $Y_t$ the stratum of $p$-rank $t$ of the $p$-rank stratification of $S$ defined by $\C$.\\

Based on (iii), one gets the following corollary:\\

\textbf{Corollary 2} {\it Let $S$ be a locally Noetherian $\F_p-$scheme and $\C$ be an $F-$crystal over $S$, then the $p$-rank strata of $S$ defined by $\C$ are pure in $S$.}\\

{\bf Proof of Corollary 2:} Let $Y_t$ be a stratum of the $p-$rank stratification of $S$ as above. A point $s\in S$ belongs to $Y_t$ if and only if $(t,0)$ is a break point of the Newton polygon of $\C_s$. Using the same notation as in the statement of Theorem 5, this shows that $Y_t=S_{(t,0)}$. By Theorem 5, $Y_t$ is affine in $S$. $\blacksquare$\\

Since purity implies weak purity, Corollary 2 implies the following theorem:\\

\textbf{Theorem 6 (Th. Zink, [Zi01, (Proposition 5)])} {\it Let $S$ be a locally Noetherian $\F_p-$scheme and $\C$ be an $F-$crystal over $S$, then the $p$-rank strata of $S$ defined by $\C$ are weakly pure in $S$.}\\

\noindent Fourth application of Theorem 5:\\

Let $R$ be an $\F_p-$algebra and let $A$ be an $n\times n$ matrix with coefficients in $R$, i.e., $A\in M_{n\times n}(R).$ Let $\underline x=\left[\begin{array}{cccc}  
x_1\\  
 x_2\\
\cdots\\
x_n\\
    \end{array}\right]$ be an $n$-dimensional vector whose entries are variables $x_1, x_2,\cdots,x_n$ and let $\underline b$ be a constant $n-$dimensional vector with coefficients in $R$, i.e., $\underline b \in R^n$. Let $\underline x^{[p]}=\left[\begin{array}{cccc}  
x_1^p\\  
 x_2^p\\
\cdots\\
x_n^p\\
    \end{array}\right]$. Consider the following Artin-Schreier equation:\\

$$\underline x=A\underline x^{[p]}+\underline b.$$\\

Let $S$ be a $R$-scheme. Define 
\begin{center}
$\phi:~S \to \N$ \end{center}
 \begin{center} $s \mapsto d= \dim_{\F_p}(\{ \underline x\in \overline{k(s)}^n|\underline x= A_s \underline{x}^{[p]}+ {\underline b_s}\})$,
\end{center}
where $\overline {k(s)}$ is the algebraic closure of the residue field $k(s)$ of $s\in S$, $A_s$ is the canonical image of $A$ in $M_{n\times n}(\overline {k(s)})$ and ${\underline b_s}$ is the canonical image of $\underline b$ in $\overline {k(s)}^n$.\\

Let $S$ be an $R$-scheme. For any $d\in \N$, let $Y_d$ be the reduced, locally closed subscheme of $S$ formed by those points $s\in S$ with the property $\phi(s)=d$. We call $Y_d$ the stratum of the Artin-Schreier stratification of $S$ defined by the equation $\underline x=A\underline x^{[p]}+\underline b$.\\

By the first proof of Theorem 2.4.1 (b) in [Va13], the Artin-Schreier stratification of $S$ defined by the equation $\underline x=A\underline x^{[p]}+\underline b$ is equivalent to the Artin-Schreier stratification of $S$ defined by the equation $\underline x=\tilde A\underline x^{[p]}$ for some $\tilde A\in M_{n\times n}(R)$. Therefore from now on, we will always assume $\underline b=\underline 0$ in the Artin-Schreier stratification.\\

\textbf{Corollary 3} {\it Let $S=\Spec~R$ be a locally Noetherian affine $\F_p$-scheme. Let $\underline x,~\underline x^{[p]}$ and $A$ be defined as above. Then each Artin-Schreier stratum of $S$ defined by the equation $\underline x=A\underline x^{[p]}$ is pure in $S$.}\\

{\bf Proof of Corollary 3:} Let $S_0$ be a stratum of the Artin-Schreier stratification of $S$ defined by the equation $\underline x=A\underline x^{[p]}$. By replacing $S$ by the schematic closure $\overline{S_0}$ of $S_0$ in $S$, we can assume that $S_0$ is open dense in $S$. We can further assume that the $\F_p$-scheme $S$ is reduced. Consider the pullback of $S_0$ in the following commutative diagram:
\begin{center}
$\xymatrixcolsep{5pc}\xymatrix{S_0^{\perf}\ar@{_{(}->}[d]^{\mbox{open}}\ar[r]^{\mbox{integral}}& S_0\ar@{_{(}->}[d]^{\mbox{open}}
\\S^{\perf}=\Spec~ R^{\perf}\ar@{->}[r]^{\mbox{integral}}&S=\Spec~ R.}$
\end{center}

Since the morphism $S_0^{\perf}\to S_0$ is integral and $S_0$ is an open subscheme of $S$, to prove that $S_0$ is affine it suffices to prove that $S_0^{\perf}$ is affine, cf. [Va06, (Lemma 2.9.2)]. Therefore, we can assume $R=R^{\perf}$, i.e., $R$ is a reduced, Noetherian perfect $\F_p-$algebra.\\ 

Now consider the $F$-crystal $\C=(W(R)^{2n},F)$ over $R$, where $F=g
\left[\begin{array}{cc}  
I_n& 0_n\\  
 0_n &pI_n\\
    \end{array}\right] \sigma_{W(R)}$. Here $g$ is a fixed invertible matrix in $GL_{2n}(W(R))$ lifting $\overline g=\left[\begin{array}{cc}  
A & I_n\\  
 I_n & 0_n\\
    \end{array}\right]\in GL_{2n}(R)$, $I_n$ is the $n\times n$ identity matrix, $0_n$ is the $n\times n$ zero matrix and $\sigma_{W(R)}$ is the Frobenius endomorphism of $W(R)$.\\

{\bf Claim:} {\it For any $d\in \N$, the stratum of $p$-rank $d$ of the $p$-rank stratification of $S$ defined by $\C$ is the same as the stratum $Y_d$ (using the same notation as in the definition of Artin-Schreier stratification) of the Artin-Schreier stratification of $S$ defined by the equation $\underline x=A\underline x^{[p]}$.}\\

{\bf Proof of Claim}: Let $s\in S=\Spec~R$, $\overline {k(s)}$ be the algebraic closure of its residue field, $\sigma_{W(\overline{k(s)})}$ be the Frobenius endomorphism of $W(\overline{k(s)})$ and $F_s=F\otimes \sigma_{W(\overline{k(s)})}=g_s
\left[\begin{array}{cc}  
I_n& 0_n\\  
 0_n &pI_n\\
    \end{array}\right] \sigma_{W(\overline{k(s)})}$, where $g_s$ is the canonical image of $g$ in $GL_{2n}(W(\overline{k(s)}))$. Consider the pullback $\C_{s}=(W(\overline{k(s)})^{2n}, F_s)$ of $\C$ to $\Spec~ \overline {k(s)}$:
\begin{center}
$\xymatrix{\C_s\ar[d]\ar[r]& \C \ar[d]\\
 \Spec~ \overline{k(s)}\ar@{->}[r]&S=\Spec ~R.}$
\end{center}

We investigate the $\F_p$-Vector space $V=(\{\underline z\in \overline \C_s|\overline {F_s}(\underline z)=\underline z\})$, where $\overline {F_s}=\overline {g_s}
\left[\begin{array}{cc}  
I_n& 0_n\\  
 0_n & 0_n\\
    \end{array}\right] \sigma_{\overline{k(s)}}=\left[\begin{array}{cc}  
A_s & 0_n\\  
 I_n & 0_n\\
    \end{array}\right] \sigma_{\overline{k(s)}} $ is the reduction modulo $p$ of $F_s$, where $A_s$ is the canonical image of $A$ in $M_{n\times n}(\overline{k(s)})$ and $\overline{g_s}$ is $g_s$ mod $p$, and $\overline{\C_s}=(\overline {k(s)}^{2n}, \overline{F_s})$ is the reduction modulo $p$ of $\C_s$. Suppose $\underline z=\left[\begin{array}{cccc}  
z_1\\  
 z_2\\
\cdots\\
z_{2n}\\
    \end{array}\right]\in \overline {k(s)}^{2n}$. We have $\underline z\in V$ if and only if $\left[\begin{array}{cc}  
A_s & 0_n\\  
 I_n & 0_n\\
    \end{array}\right]\left[\begin{array}{cccc}  
z_1^p\\  
 z_2^p\\
\cdots\\
z_{2n}^p\\
    \end{array}\right]=\left[\begin{array}{cccc}  
z_1\\  
 z_2\\
\cdots\\
z_{2n}\\
    \end{array}\right]$ if and only if $A_s\left[\begin{array}{cccc}  
z_1^p\\  
 z_2^p\\
\cdots\\
z_n^p\\
    \end{array}\right]=\left[\begin{array}{cccc}  
z_1\\  
 z_2\\
\cdots\\
z_n\\
    \end{array}\right]$ and $z_{n+1}=z_1^p,z_{n+2}=z_2^p,\cdots,z_{2n}=z_n^p$. Therefore $$\dim_{\F_p}V=\dim_{\F_p}(\{ \underline x\in \overline{k(s)}^n| \underline x=A_s \underline{x}^{[p]}\}):=d,$$which means the stratum of $p$-rank $d$ in the $p$-rank stratification of $S$ defined by $\C$ is the same as the stratum $Y_d$ in the Artin-Schreier stratification of $S$ defined by the equation $\underline x=A\underline x^{[p]}$. Since $s\in S$ is arbitrary chosen and thus $d$ is arbitrary, Claim is proved.\\

By Corollary 2, each Artin-Schreier stratum of $S$ defined by the equation $\underline x=A\underline x^{[p]}$ is pure in $S$. This ends the proof of Corollary 3. $\blacksquare$\\

Remark: Deligne and Vasiu also obtained the results of Corollary 2 and 3 using different methods, cf. [De11] and [Va14].

\newpage

\noindent\textbf{\Huge{References}}\\

\noindent[Be74] Pierre Berthelot: {\it Cohomologie cristalline des sch\'emas de caract\'eristique $p > 0$}, Lecture Notes in Math., Vol. 407, Springer-Verlag, Berlin-New York, 1974.\\

\noindent[BM90] Pierre Berthelot and William Messing: {\it Th\'eorie de Dieudonn\'e cristalline III. Th\'eor\`emes d'\'equivalence et de pleine fid\'elit\'e}, The Grothendieck Festschrift, Vol. I, 173–247, Progr. Math., Vol. 86, Birkh\"{a}user Boston, Boston, MA, 1990.\\

\noindent[De11] Pierre Deligne: An unpublished proof, notes to A. Vasiu, 2011.\\

\noindent[De72] Michel Demazure: {\it Lectures on p-divisible groups}, Lecture Notes in Math., vol. 302, Springer-Verlag, 1972.\\

\noindent[Ei07] David Eisenbud: {\it Commutative algebra with a view toward algebraic geometry}, Grad. Texts in Math., vol. 150, Springer-Verlag, 2007.\\ 

\noindent[Gr61] Alexander Grothendieck: {\it \'El\'ements de g\'eom\'etrie alg\'ebrique. II. \'Etude globale \'el\'ementaire de quelques classes de morphisms}, Inst. Hautes \'Etudes Sci. Publ. Math., vol. 11, 1961.\\

\noindent[Gr64] Alexander Grothendieck: {\it \'El\'ements de g\'eom\'etrie alg\'ebrique. IV. \'Etude locale des sch\'emas et desmorphismes de sch\'ema}, Inst. Hautes \'Etudes Sci. Publ. Math., vol. 20, 1964, vol. 24, 1965, vol. 28, 1996, and vol. 32, 1967.\\

\noindent[Ha08] Michiel Hazewinkel: {\it Witt vectors}, part 1, manuscript, 148 pages, 2008.\\

\noindent[Hi80] Matsumura Hideyuki: {\it Commutative algebra, second ed.}, Benjamin/Cummings, Reading, MA, 1980.\\

\noindent[HT01] Michael Harris and Richard Taylor: {\it On the geometry and cohomology of some simple Shimura varieties}, {\em Annals of Math. Studies} 151, PUP 2001.\\

\noindent[Il94] Luc Illusie: {\it Crystalline cohomology. Motives (Seattle, WA, 1991), 43–70}, {\em Proc. Sympos. Pure Math.}, 55, Part 1, Amer. Math. Soc., Providence, RI, 1994.\\ 

\noindent[JO00] Aise Johan de Jong and Frans Oort: {\it Purity of the stratification by Newton polygons}, {\em J. Amer. Math. Soc.} 13 (2000), no. 1, 209--241.\\

\noindent[Ka79] Nicholas Michael Katz: {\it Slope filtration of F-crystals}, {\em Journ\'ees de G\'eom\'etrie Alg\'ebrique de Rennes} (Rennes, 1978), Vol. I, Ast\'erisque, No. 63 (1979), 113--163.\\

\noindent[NVW10] Marc-Hubert Nicole, Adrian Vasiu and Torsten Wedhorn: {\it Purity of level m stratifications}, {\em Ann. Sci. \'Ecole Norm. Sup. }(4) 43 (2010), no. 6, 925–955.\\

\noindent[Se79] Jean-Pierre Serre: {\it Local fields}, Grad. Texts in Math., vol. 67, Springer-Verlag, 1979.\\

\noindent[Va06] Adrian Vasiu: {\it Crystalline boundedness principle}, {\em Ann. Sci. \'Ecole Norm. Sup.}  39 (2006), no. 2, 245--300.\\

\noindent[Va13] Adrian Vasiu: {\it A motivic conjecture of Milne}, {\em J. Reine Agew. Math. (Crelle)} 685 (2013), 181–247.\\

\noindent[Va14] Adrian Vasiu: {\it Talk on Conference on Arithmetic Algebraic Geometry on the occasion of Gerd Faltings' 60th birthday}, Max Planck Institute for Mathematics, 2014.\\

\noindent[Ya10] Yanhong Yang: {\it An improvement of de Jong--Oort's purity theorem}, preprint, 2010, arXiv:1004.3090.\\

\noindent[Zi01] Thomas Zink: {\it On the slope filtration}, {\em Duke Math. J.}  109 (2001), no. 1, 79--95.\\

\end{document}